\documentclass[11pt]{article}%
\usepackage{graphicx}
\usepackage{amsmath}
\usepackage{enumerate}%
\usepackage{amsfonts}%
\usepackage{amssymb}
\newtheorem{theorem}{Theorem}[section]
\newtheorem{acknowledgement}{Acknowledgement}

\newtheorem{corollary}{Corollary}[section]

\newtheorem{lemma}{Lemma}[section]

\newtheorem{proposition}{Proposition}[section]

\newenvironment{proof}[1][Proof]{\textbf{#1.} }{\ \rule{0.5em}{0.5em}}

\begin{document}

\title{Flow properties of differential equations driven by fractional Brownian motion}
\author{L. Decreusefond and D. Nualart}
\date{}
\maketitle

\begin{abstract}
We prove that solutions of stochastic differential equations driven by
fractional Brownian motion for $H>1/2$ define flows of homeomorphisms on $\mathbb{R}^{d}$.
$\smallskip$\newline \textbf{AMS Subject Classification: }60H05, 60H07

\end{abstract}

\section{Introduction}

Suppose that $B^{H}=\{B_{t}^{H},t\geq0\}$ is an $m$-dimensional
\textit{fractional Brownian motion} with Hurst parameter $H\in(0,1)$, defined
in a complete probability space $(\Omega,\mathcal{F},P)$. That is, the
components $B^{H,j}$, $j=1,\ldots,m$, are independent zero mean Gaussian
processes with the covariance function%
\begin{equation}
R_{H}(t,s)=\frac{1}{2}\left(  t^{2H}+s^{2s}-\left|  t-s\right|  ^{2H}\right)
.\label{f1}%
\end{equation}
For $H=\frac{1}{2}$, the process $B^{H}$ is an $m$-dimensional ordinary
Brownian motion. On the other hand, from (\ref{f1}) it follows that
\[
E(\left|  B_{t}^{H,j}-B_{s}^{H,j}\right|  ^{2})=|t-s|^{2H}.
\]
As a consequence, the processes $B^{H,j}$ have stationary increments, and for
any $\alpha<H$ we can select versions with H\"{o}lder continuous trajectories
of order $\alpha$ on a compact interval $[0,T]$. 

This process was first studied by Kolmogorov in \cite{Ko} and later by
Mandelbrot and Van Ness in \cite{Ma}, where a stochastic integral
representation in terms of an ordinary Brownian motion was established.

The fractional Brownian motion has the following \textit{self-similar}
property: For any constant $a>0$, the processes $\left\{  a^{-H}B_{at}%
^{H},t\geq0\right\}  $ and $\left\{  B_{t}^{H},t\geq0\right\}  $ have the same
distribution. For $H=\frac{1}{2}$ the process $B^{H}$ has independent
increments, but for $H\neq\frac{1}{2}$, this property is no longer true. In
particular, if $H>\frac{1}{2}$, the fractional Brownian motion has the long range dependence
property, which means that for each $j=1,\ldots,m$%
\[
\sum_{n=1}^{\infty}\mathrm{Corr}(B_{n+1}^{H,j}-B_{n}^{H,j},B_{1}^{H,j}%
)=\infty.
\]

The self-similar and long range dependence properties make the fractional
Brownian motion a convenient model for some input noises in a variety of
topics from finance to telecommunication networks, where the Markov property
is not required. This fact has motivated the recent development of the
stochastic calculus with respect to the fractional Brownian motion. We refer
to \cite{Nu1} for a survey of this topic.

In this paper we are interested in stochastic differential equations on
$\mathbb{R}^{d}$ driven by a multi-dimensional fractional Brownian motion with
Hurst parameter $H>\frac{1}{2}$, that is, equations of the form
\begin{equation}
X_{t}^{i}=X_{0}^{i}+\sum_{j=1}^{m}\int_{0}^{t}\sigma^{i,j}(s,X_{s}%
)dB_{s}^{H,j}+\int_{0}^{t}b^{i}(s,X_{s})ds,\label{f2}%
\end{equation}
$i=1,\ldots,d$. The stochastic integral appearing in (\ref{f2}) is a
path-wise Riemann-Stieltjes integral. In fact, under suitable conditions on
$\sigma$, the processes $\sigma(s,X_{s})$ and $B_{s}^{H}$ have trajectories
which are H\"{o}lder continuous of order strictly larger than $\frac{1}{2}$ and
we can use the approach introduced by Young in \cite{Yo}. A first result on
the existence and uniqueness of a solution for this kind of equations was
obtained by Lyons in \cite{Ly1}, using the notion of $p$-variation. On the
other hand, the theory of rough path analysis introduced by Lyons in \cite{Ly1}
\ (see also the monograph by Lyons and Qian   \cite{LQ}), has allowed Coutin
and Qian \cite{cq} to establish the existence of strong solutions and a
Wong-Zakai type approximation limit for the stochastic differential equations
\ of the form\ (\ref{f2}) driven by a fractional Brownian motion with
parameter $H>\frac{1}{4}$. In \cite{cq}  sufficient conditions on the vector
fields $b$ and $\sigma$ are given to ensure existence and uniqueness of the
solution of (\ref{f2}) even when vector fields do not commute.

In \cite{Za} Z\"{a}hle has introduced a generalized Stieltjes integral using
the techniques of fractional calculus. This integral is expressed in terms of
fractional derivative operators and it coincides with the Riemann-Stieltjes
integral $\int_{a}^{b}fdg$, when the functions $f$ and $g$ are H\"{o}lder
continuous of orders $\lambda$ and $\mu$, respectively and $\lambda+\mu>1$.
Using this formula for the Riemann-Stieltjes integral, Nualart and
R\u{a}\c{s}canu have obtained in \cite{NR} the existence of a unique solution
for the stochastic differential equations (\ref{f2}) under general conditions
on the coefficients.

Later on, Nualart and Saussereau \cite{NS} have studied the regularity in the
sense of Malliavin Calculus of the solution of Equation (\ref{f2}), and they
have established the absolute continuity of the law of the random variable
$X_{t}$ under some non-degeneracy conditions on the coefficient $\sigma$.

The main result of this paper is the flow and homeomorphic properties of $X$ as
a function of the initial condition $x$. Since the solution of (\ref{f2}) is
defined path-wise, ordinary (i.e., deterministic) methods are in use here.
Namely, we use the estimates found in \cite{NR} and approximate fractional
Brownian motion by a sequence of regular processes  to prove that the solution
$\{X_{rt}(x),\,0\leq r\leq t\leq T,\,x\in\mathbb{R}^{d}\}$ of
\[
X_{rt}(x)=x+\int_{r}^{t}\sigma\left(  s,X_{rs}(x)\right)  dB^{H}(s)+\int
_{r}^{t}b(s,X_{rs}(x))ds,
\]
defines a flow of ${}\mathbb{R}^{d}$-homeomorphisms.

The paper is organized as follows. Section 2 contains some preliminaries on
fractional calculus. In Section 3 we review some results on the properties of
the solution of Equation (\ref{f2}) and we establish some continuity estimates
as a function of the initial condition and the driven input, which are needed
later. finally in Section 4 we show that Equation (\ref{f2}) defines a flow of homeomorphisms.

\section{Preliminaries}

\label{sec:preliminaries}

Let $a,b\in{},$ $a<b.$ Let $f\in L^{1}\left(  a,b\right)  $ and $\alpha>0.$
The left-sided and right-sided fractional Riemann-Liouville integrals of $f$
of order $\alpha$ are defined for almost all $x\in\left(  a,b\right)  $ by
\[
I_{a+}^{\alpha}f\left(  x\right)  =\frac{1}{\Gamma\left(  \alpha\right)  }%
\int_{a}^{x}\left(  x-y\right)  ^{\alpha-1}f\left(  y\right)  dy
\]
and
\[
I_{b-}^{\alpha}f\left(  x\right)  =\frac{\left(  -1\right)  ^{-\alpha}}%
{\Gamma\left(  \alpha\right)  }\int_{x}^{b}\left(  y-x\right)  ^{\alpha
-1}f\left(  y\right)  dy,
\]
respectively, where $\left(  -1\right)  ^{-\alpha}=e^{-i\pi\alpha}$ and
$\Gamma\left(  \alpha\right)  =\int_{0}^{\infty}r^{\alpha-1}e^{-r}dr$ is the
Euler function. Let $I_{a+}^{\alpha}(L^{p})$ (resp. $I_{b-}^{\alpha}(L^{p})$)
the image of $L^{p}(a,b)$ by the operator $I_{a+}^{\alpha}$ (resp.
$I_{b-}^{\alpha}$). If $f\in I_{a+}^{\alpha}\left(  L^{p}\right)  \ $ (resp.
$f\in I_{b-}^{\alpha}\left(  L^{p}\right)  $) and $0<\alpha<1$ then the Weyl
derivative
\begin{equation}
D_{a+}^{\alpha}f\left(  x\right)  =\dfrac{1}{\Gamma\left(  1-\alpha\right)
}\left(  \dfrac{f\left(  x\right)  }{\left(  x-a\right)  ^{\alpha}}+\alpha
\int_{a}^{x}\dfrac{f\left(  x\right)  -f\left(  y\right)  }{\left(
x-y\right)  ^{\alpha+1}}dy\right)  1_{\left(  a,b\right)  }(x) \label{a2}%
\end{equation}%
\begin{equation}
\left(  \text{resp. \ \ }D_{b-}^{\alpha}f\left(  x\right)  =\dfrac{\left(
-1\right)  ^{\alpha}}{\Gamma\left(  1-\alpha\right)  }\left(  \dfrac{f\left(
x\right)  }{\left(  b-x\right)  ^{\alpha}}+\alpha\int_{x}^{b}\dfrac{f\left(
x\right)  -f\left(  y\right)  }{\left(  y-x\right)  ^{\alpha+1}}dy\right)
1_{\left(  a,b\right)  }(x)\right)  \label{a3}%
\end{equation}
is defined for almost all $x\in\left(  a,b\right)  .$

For any $0<\lambda\leq1$, denote by $C^{\lambda} (  0,T )  $ the
space of $\lambda$-H\"{o}lder continuous functions $f:[0,T]\rightarrow
\mathbb{R}$, equipped with the norm $\left\|  f\right\|  _{\infty
}+\left\|  f\right\|  _{\lambda}$, where%
\[
\left\|  f\right\|  _{\lambda}=\sup_{0\leq s<t\leq T}\dfrac{\left|  f\left(
t\right)  -f\left(  s\right)  \right|  }{\left(  t-s\right)  ^{\lambda}}.
\]

Recall from \cite{samko93} that we have:

\begin{itemize}
\item If $\alpha<\dfrac{1}{p}$ and $q=\dfrac{p}{1-\alpha p}$ then $
I_{a+}^{\alpha}\left(  L^{p}\right)  =I_{b-}^{\alpha}\left(  L^{p} \right)
\subset L^{q}\left(  a,b\right)  . $

\item If $\alpha>\dfrac{1}{p}$ then$I_{a+}^{\alpha}\left(  L^{p}\right)
\,\cup\,I_{b-}^{\alpha}\left(  L^{p}\right)  \subset C^{\alpha-\frac{1}{p}%
}\left(  a,b\right)  . $
\end{itemize}

The linear spaces $I_{a+}^{\alpha}\left(  L^{p}\right)  $ are Banach spaces
with respect to the norms
\[
\left\|  f\right\|  _{I_{a+}^{\alpha}\left(  L^{p}\right)  }=\left\|
f\right\|  _{L^{p}}+\left\|  D_{a+}^{\alpha}f\right\|  _{L^{p}}\thicksim
\left\|  D_{a+}^{\alpha}f\right\|  _{L^{p}},
\]
and the same is true for $I_{b-}^{\alpha}\left(  L^{p}\right)  $.

Suppose that $f\in C^{\lambda}(a,b)$ and $g\in C^{\mu}(a,b)$ with $\lambda
+\mu>1$. Then, from the classical paper by Young \cite{Yo}, the
Riemann-Stieltjes integral $\int_{a}^{b}fdg$ exists. The following proposition
can be regarded as a fractional integration by parts formula, and provides an
explicit expression for the integral $\int_{a}^{b}fdg$ in terms of fractional
derivatives (see \cite{Za}).

\begin{proposition}
\label{p1} Suppose that $f\in C^{\lambda}(a,b)$ and $g\in C^{\mu}(a,b)$ with
$\lambda+\mu>1$. Let ${\lambda}>\alpha$ and $\mu>1-\alpha$. Then the Riemann
Stieltjes integral $\int_{a}^{b}fdg$ exists and it can be expressed as%
\begin{equation}
\int_{a}^{b}fdg=(-1)^{\alpha}\int_{a}^{b}D_{a+}^{\alpha}f\left(  t\right)
D_{b-}^{1-\alpha}g_{b-}\left(  t\right)  dt, \label{1.8}%
\end{equation}
where $g_{b-}\left(  t\right)  =g\left(  t\right)  -g\left(  b\right)  $.
\end{proposition}

In \cite{Za} Z\"{a}hle introduced a generalized Stieltjes integral of
$f$ with respect to $g$ defined by the right-hand side of (\ref{1.8}),
assuming that $f$ and $g$ are functions such that $\ $ $g(b-)$ exists, $f\in
I_{a+}^{\alpha}\left(  L^{p}\right)  $ and $g_{b-}\in I_{b-}^{1-\alpha}\left(
L^{q}\right)  $ for some $p,q\geq1,$ $1/p+1/q\leq1,$ $0<\alpha<1$.

\medskip
Let $\ \alpha<\frac{1}{2}$ and $d\in\mathbb{N}^{\ast}.$ Denote by
$W_{0}^{\alpha,\infty}(0,T;\mathbb{R}^{d})$ the space of measurable functions
$f:[0,T]\rightarrow\mathbb{R}^{d}$ such that%

\[
\left\|  f\right\|  _{\alpha,\infty}:=\sup_{t\in\lbrack0,T]}\left(
|f(t)|+\int_{0}^{t}\dfrac{\left|  f\left(  t\right)  -f\left(  s\right)
\right|  }{\left(  t-s\right)  ^{\alpha+1}}ds\right)  <\infty.
\]
We have, for all $0<\varepsilon<\alpha$%
\[
C^{\alpha+\varepsilon}(0,T;\mathbb{R}^{d})\subset W_{0}^{\alpha,\infty
}(0,T;\mathbb{R}^{d})\subset C^{\alpha-\varepsilon}(0,T;\mathbb{R}^{d}).
\]
Denote by $W_{T}^{1-\alpha,\infty}(0,T;\mathbb{R}^{m})$ the space of
measurable functions $g:[0,T]\rightarrow\mathbb{R}^{m}$ such that%
\[
\left\|  g\right\|  _{1-\alpha,\infty,T}:=\sup_{0<s<t<T}\left(
\frac{|g(t)-g(s)|}{(t-s)^{1-\alpha}}+\int_{s}^{t}\frac{|g(y)-g(s)|}%
{(y-s)^{2-\alpha}}dy\right)  <\infty.
\]
Clearly, for all $\varepsilon>0$ we have
\begin{equation}
C^{1-\alpha+\varepsilon}\left(  0,T;\mathbb{R}^{m}\right)  \subset W_{T}^{1-\alpha,\infty
}(0,T;\mathbb{R}^{m})\subset C^{1-\alpha}\left(  0,T;\mathbb{R}^{m}\right)  .\label{eq:2}%
\end{equation}
Moreover, if $g\ $belongs to $W_{T}^{1-\alpha,\infty}(0,T;\mathbb{R}^{m}),$ its restriction
to $(0,t)$ belongs to $I_{t-}^{1-\alpha}(L^{\infty}(0,t;\mathbb{R}^{m}))$ for all $t$ and
\begin{align}
\Lambda_{\alpha}(g) &  :=\frac{1}{\Gamma(1-\alpha)}\sup_{0<s<t<T}\left|
\left(  D_{t-}^{1-\alpha}g_{t-}\right)  (s)\right|  \nonumber\\
&  \leq\frac{1}{\Gamma(1-\alpha)\Gamma(\alpha)}\left\|  g\right\|
_{1-\alpha,\infty,T}<\infty.\label{eq:3}%
\end{align}
The integral $\int_{0}^{t}fdg$ can be defined for all $t\in\lbrack0,T]$ if $g$
belongs to $W_{T}^{1-\alpha,\infty}(0,T)$ and $f$ satisfies%
\[
\left\|  f\right\|  _{\alpha,1}:=\int_{0}^{T}\frac{|f(s)|}{s^{\alpha}}%
ds+\int_{0}^{T}\int_{0}^{s}\frac{|f(s)-f(y)|}{(s-y)^{\alpha+1}}dy\ ds<\infty.
\]
Furthermore the following estimate holds
\[
\left|  \int_{0}^{T}fdg\right|  \leq\Lambda_{\alpha}(g)\left\|  f\right\|
_{\alpha,1}.
\]

\section{Stochastic differential equations driven by a fBm}

We are going to consider first the case of a deterministic equation. Let
$0<\alpha<\frac{1}{2}$ be fixed. Let $g\in W_{T}^{1-\alpha,\infty}%
(0,T;{}\mathbb{R}^{m})$. Consider the deterministic differential equation on
${}\mathbb{R}^{d}$
\begin{equation}
\xi_{t}^{i}=x_{0}^{i}+\int_{0}^{t}b^{i}(s,\xi_{s})ds+\sum_{j=1}^{m}\int
_{0}^{t}\sigma^{i,j}\left(  s,\xi_{s}\right)  dg_{s}^{j},\,\;t\in\left[
0,T\right]  , \label{2.1}%
\end{equation}
$i=1,...,d$, where $x_{0}\in{}{}\mathbb{R}^{d}$, and the coefficients
$\sigma^{i,j},b^{i}:\left[  0,T\right]  \times{}{}\mathbb{R}^{d}%
\mathbb{\rightarrow R}{}$ are measurable functions. Set $\sigma=\left(
\sigma^{i,j}\right)  _{d\times m},\,\;\,b=\left(  b^{i}\right)  _{d\times1}$
and for a matrix $A=\left(  a^{i,j}\right)  _{d\times m}$ and a vector
$y=\left(  y^{i}\right)  _{d\times1}$ denote $\left|  A\right|  ^{2}%
=\sum_{i,j}\left|  a^{i,j}\right|  ^{2}$ and $\left|  y\right|  ^{2}=\sum
_{i}\left|  y^{i}\right|  ^{2}.$

Let us consider the following assumptions on the coefficients.

\begin{enumerate}
\item[(H1)] $\sigma(t,x)$ is differentiable in $x$, and \ there exist some
constants $\ $ $0<\beta,\delta\leq1$, $M_{1},M_{2},M_{3}>0$ such that the
following properties hold:%
\[
\left\{
\begin{array}
[c]{l}%
i)\quad\text{Lipschitz continuity}\\
|\sigma(t,x)-\sigma(t,y)|\leq M_{1}|x-y|,\quad\forall\,x\in{}\mathbb{R}%
^{d},\,\forall\,t\in\left[  0,T\right]  \smallskip\\
ii)\quad\text{H\"{o}lder continuity}\\
|\partial_{x_{i}}\sigma(t,x)-\partial_{x_{i}}\sigma(t,y)|\leq M_{2}%
|x-y|^{\delta}\smallskip,\\
 \quad\quad\quad\quad\quad\quad\quad\quad\;\;\forall\left|  x\right|
,\left|  y\right|  \in{}\mathbb{R}^{d},\,\forall t\in\left[  0,T\right], i=1,\dots,d
,\smallskip\\
iii)\quad\text{H\"{o}lder continuity in time}\\
\left|  \sigma(t,x)-\sigma(s,x)\right|  +|\partial_{x_{i}}\sigma
(t,x)-\partial_{x_{i}}\sigma(s,x)|\leq M_{3}|t-s|^{\beta}\\
\quad\quad\quad\quad\quad\quad\quad\quad\quad\quad\quad\quad\forall\,x\in
{}\mathbb{R}^{d},\,\forall\,t,s\in\left[  0,T\right].
\end{array}
\right.
\]

\item[(H2)] There exists constants $L_{1}, L_{2}>0$ such that the following
properties hold:%
\[
\left\{
\begin{array}
[c]{ll}%
i) & \quad\text{Local Lipschitz continuity}\\
& \left|  b(t,x)-b(t,y)\right|  \leq L_{1}|x-y|,\;\forall\left|  x\right|
,\left|  y\right|  \in{}\mathbb{R}^{d},\,\forall t\in\left[  0,T\right]
,\smallskip\\
ii) & \quad\text{Linear growth}\\
& \left|  b(t,x)\right|  \,\leq\,L_{2}(1+|x|),\;\forall x\in{}\mathbb{R}%
^{d},\,\forall t\in\left[  0,T\right]  \smallskip.
\end{array}
\right.
\]
Set%
\[
\alpha_{0}=\min\left\{  \frac{1}{2},\beta,\frac{\delta}{1+\delta}\right\}  .
\]
The following existence and uniqueness result has been proved in
\cite{NR}.
\end{enumerate}

\begin{theorem}
Suppose that the coefficients $\sigma(t,x)$ and $b(t,x)$ satisfy assumptions
(H1) and (H2). Then, if $\alpha<\alpha_{0}$ there exists a unique solution of
Equation (\ref{2.1}) in the space $C^{1-\alpha}\left(  0,T;\ \mathbb{R}%
^{d}\right)  .$
\end{theorem}

Actually, these conditions can be slightly relaxed. For instance, the
H\"{o}lder continuity of the partial derivatives of $\sigma$ and the
H\"{o}lder continuity of the coefficient $b$ may hold only locally (see \cite{NR} for the details).

We now state two theorems which are consequences of the estimates found in
\cite{NR}.

For any $\lambda\geq0$ we introduce the equivalent norm \ in \ the space
$W_{0}^{\alpha,\infty}\left(  0,T;\mathbb{R}^{d}\right)  $ defined by%
\[
\left\|  f\right\|  _{\alpha,\lambda}=\sup_{t\in\lbrack0,T]}\ e^{-\lambda
t}\left(  \left|  f(t)\right|  +\int_{0}^{t}\frac{|f(t)-f(s)|}{(t-s)^{\alpha
+1}}\ \ ds\right)  .
\]

\begin{theorem}
\label{thm:continuite_condition_initiale} Let us denote $\xi_{t}(x_{0})$ the
solution of (\ref{2.1}) at time $t$ with initial condition $x_{0}$. Fix $R>1$.
Then there exists a constant $C$ such that for any $x_{0}$ and $x_{1}$ in the
ball $B(0,\,R)=\{x\in{}\mathbb{R}^{d},\ |x|\leq R\}$ and for any $\lambda>$
$\left[  R\exp\left(  C\left(  1+\Lambda_{\alpha}(g)\right)  ^{\frac{1}%
{1-2\alpha}}\right)  \right]  ^{\frac{1}{1-2\alpha}}$ we have
\[
\lVert\xi(x_{0})-\xi(x_{1})\lVert_{\alpha,\,\lambda}\leq\left(  1-R\exp\left(
C\left(  1+\Lambda_{\alpha}(g)\right)  ^{\frac{1}{1-2\alpha}}\right)
\lambda^{2\alpha-1}\right)  ^{-1}\ |x_{0}-x_{1}|.
\]
\end{theorem}

\begin{proof}
It is proved in \cite{NR} that there exists a constant $C_{1}$ such
that if $\lambda_{0}=C_{1}\left(  1+\Lambda_{\alpha}(g)\right)  ^{\frac{1}%
{1-2\alpha}}$ then for any initial condition $x_{0}$ in the ball of radius $R$
we have
\begin{equation}
\left\|  \xi(x_{0})\right\|  _{\alpha,\lambda_{0}}\leq2\left(  1+\left|
x_{0}\right|  \right)  \leq4R. \label{f5}%
\end{equation}
Given a function $f\in W_{0}^{\alpha,\infty}(0,T;\mathbb{R}^{d})$ we define as
in \cite{NR}
\begin{align*}
F_{t}^{(b)}(f)  &  =\int_{0}^{t}b(s,\,f(s))ds,\\
G_{t}^{(\sigma)}(g,f)  &  =\int_{0}^{t}\sigma(s,\,f(s))dg(s),\\
\Delta(f)  &  =\sup_{r\in\lbrack0,T]}\int_{0}^{r}\frac{|f(r)-f(s)|^{\delta}%
}{(r-s)^{\alpha+1}}ds.
\end{align*}
If $f,h\in W_{0}^{\alpha,\infty}(0,T;\mathbb{R}^{d})$ (see \cite{NR})
there exist constants $C_{2}$ and $C_{3}$ such that%
\begin{align}
\left\|  F^{(b)}(f)-F^{(b)}(h)\right\|  _{\alpha,\lambda}  &  \leq
C_{2}\lambda^{\alpha-1}\left\|  f-h\right\|  _{\alpha,\lambda},\label{f7}\\
\left\|  G^{\left(  \sigma\right)  }\left(  g,f\right)  -G^{\left(
\sigma\right)  }\left(  g,h\right)  \right\|  _{\alpha,\lambda}  &  \leq
C_{3}\Lambda_{\alpha}(g)\lambda^{2\alpha-1}   \notag \\
   & \times \left(  1+\Delta\left(  f\right)
+\Delta\left(  h\right)  \right)  \left\|  f-h\right\|  _{\alpha,\lambda
},\label{f8}\\
\left\|  G^{\left(  \sigma\right)  }\left(  g,f\right)  \right\|
_{\alpha,\lambda}  &  \leq C_{4}\Lambda_{\alpha}(g)\lambda^{2\alpha-1}\left(
1+\left\|  f\right\|  _{\alpha,\lambda}\right)  \label{f8a}%
\end{align}
Also, if $f\in W_{0}^{\alpha,\infty}(0,T;\mathbb{R}^{d})$ and $h$ is a bounded
measurable function, then%
\begin{align}
\left\|  F^{(b)}(h)\right\|  _{1-\alpha}  &  \leq C_{5}(1+\left\|  h\right\|
_{\infty}),\label{f3}\\
\left\|  G^{\left(  \sigma\right)  }\left(  g,f\right)  \right\|  _{1-\alpha}
&  \leq C_{6}\Lambda_{\alpha}(g)\left(  1+\left\|  f\right\|  _{\alpha,\infty
}\right)  . \label{f4}%
\end{align}
We have the following estimate%
\begin{align}
\Delta(\xi(x_{0}))  &  =\sup_{r\in\lbrack0,T]}\int_{0}^{r}\frac{|\xi_{r}%
(x_{0})-\xi_{s}(x_{0})|^{\delta}}{(r-s)^{\alpha+1}}\ ds\nonumber\\
&  \leq\dfrac{T^{\delta-\alpha\left(  1+\delta\right)  }}{\delta-\alpha\left(
1+\delta\right)  }\left\|  \xi(x_{0})\right\|  _{1-\alpha}, \label{f9}%
\end{align}
and using (\ref{f3}), (\ref{f4}), and (\ref{f5}) we obtain
\begin{align*}
\left\|  \xi(x_{0})\right\|  _{1-\alpha}  &  \leq\left|  x_{0}\right|
+\left\|  F^{\left(  b\right)  }\left(  \xi(x_{0})\right)  \right\|
_{1-\alpha}+\left\|  G^{\left(  \sigma\right)  }\left(  \xi(x_{0})\right)
\right\|  _{1-\alpha}\\
&  \leq\left|  x_{0}\right|  +C_{5}\left(  1+\left\|  \xi(x_{0})\right\|
_{\infty}\right)  +\Lambda_{\alpha}(g)C_{6}\left(  1+\left\|  \xi
(x_{0})\right\|  _{\alpha,\infty}\right) \\
&  \leq C_{7}e^{\lambda_{0}T}(1+|x_{0}|)(1+\Lambda_{\alpha}(g))\\
&  \leq2C_{7}e^{\lambda_{0}T}R(1+\Lambda_{\alpha}(g))\\
&  =2C_{7}\exp\left(  TC_{1}\left(  1+\Lambda_{\alpha}(g)\right)
^{\frac{1}{1-2\alpha}}\right)  R(1+\Lambda_{\alpha}(g)).
\end{align*}
Hence, from (\ref{f9}) we get%
\begin{equation}
\Delta(\xi(x_{0}))\leq C_{8}\exp\left(  TC_{1}\left(  1+\Lambda_{\alpha
}(g)\right)  ^{\frac{1}{1-2\alpha}}\right)  R(1+\Lambda_{\alpha}(g)).
\label{f10}%
\end{equation}
Using (\ref{f7}), (\ref{f8}), and (\ref{f10}) we obtain for $x_{0}$ and
$x_{1}$ in the ball of radius $R$
\begin{align*}
\lVert\xi(x_{0})-\xi(x_{1})\rVert_{\alpha,\,\lambda}  &  \leq|x_{0}%
-x_{1}|+\lVert F^{(b)}(\xi(x_{0}))-F^{(b)}(\xi(x_{1}))\rVert_{\alpha
,\,\lambda}\\
&  +\left\|  G^{(\sigma)}(g,\xi(x_{0}))-G^{(\sigma)}(g,\xi(x_{1}))\right\|
_{\alpha,\lambda}\\
&  \leq|x_{0}-x_{1}|+C_{1}\lambda^{\alpha-1}\left\|  \xi(x_{0})-\xi
(x_{1})\right\|  _{\alpha,\lambda}\\
&  +C_{2}\Lambda_{\alpha}(g)\lambda^{2\alpha-1}\left(  1+\Delta\left(
\xi(x_{0})\right)  +\Delta\left(  \xi(x_{1})\right)  \right)  \left\|
\xi(x_{0})-\xi(x_{1})\right\|  _{\alpha,\lambda}\\
&  \leq|x_{0}-x_{1}|+C_{1}\lambda^{\alpha-1}\left\|  \xi(x_{0})-\xi
(x_{1})\right\|  _{\alpha,\lambda}\\
&  +C_{2}\Lambda_{\alpha}(g)\lambda^{2\alpha-1}\left(  1+2C_{8}\exp\left(
TC_{1}\left(  1+\Lambda_{\alpha}(g)\right)  ^{\frac{1}{1-2\alpha}}\right)
R(1+\Lambda_{\alpha}(g))\right) \\
&  \times\left\|  \xi(x_{0})-\xi(x_{1})\right\|  _{\alpha,\lambda}.
\end{align*}
As a consequence,%
\[
\lVert\xi(x_{0})-\xi(x_{1})\rVert_{\alpha,\,\lambda}\leq|x_{0}-x_{1}%
|+\lambda^{2\alpha-1}\left(  \exp\left(  C\left(  1+\Lambda_{\alpha
}(g)\right)  ^{\frac{1}{1-2\alpha}}\right)  R\right)  \left\|  \xi(x_{0}%
)-\xi(x_{1})\right\|  _{\alpha,\lambda}%
\]
for some constant $C$, which implies the result.
\end{proof}

\begin{theorem}
\label{thm:continuite} The map
\begin{align*}
\xi\ :\ W_{T}^{1-\alpha,\infty}(0,\,T;\mathbb{R}^{m})  &  \longrightarrow
W_{0}^{\alpha,\infty}(0,\,T;\mathbb{R}^{d})\\
g  &  \longmapsto\xi,
\end{align*}
where $\xi$ is the solution of (\ref{2.1}) with $x_{0}\in B(0,\,R)$, is
continuous. Namely, for $\lambda>0$ large enough we have
\[
\lVert\xi(g)-\xi(h)\rVert_{\alpha,\,\lambda}\leq\,\frac{C_{1}\lambda
^{2\alpha-1}\lVert\xi(g)\rVert_{\alpha,\,\lambda}}{1-C_{2}\lambda^{2\alpha
-1}\ (1+\Lambda_{\alpha}(h))}\Lambda_{\alpha}(g-h).
\]
\end{theorem}

\begin{proof}
With the previous notations, we have
\[
\xi(g)=x_{0}+F^{(b)}(\xi(g))+G^{(\sigma)}(g,\xi(g))\text{ }%
\]
and%
\[
\xi(h)=x_{0}+F^{(b)}(\xi(h))+G^{(\sigma)}(h,\,\xi(h)).
\]
Thus,
\begin{align*}
\xi(g)-\xi(h)  &  =F^{(b)}(\xi(g))-F^{(b)}(\xi(h))\\
&  +G^{(\sigma)}(g-h,\,\xi(g))+\left(  G^{(\sigma)}(h,\xi(g))-G^{(\sigma
)}(h,\,\xi(h))\right)  .
\end{align*}
According to (\ref{f8a}), for $\lambda$ sufficiently large,
\begin{align*}
\lVert\xi(g)-\xi(h)\rVert_{\alpha,\,\lambda}  &  \leq C_{2}\lambda^{\alpha
-1}\left\|  \xi(g)-\xi(h)\right\|  _{\alpha,\lambda}\\
&  +C_{4}\lambda^{2\alpha-1}(1+\Lambda_{\alpha}(g-h))\lVert\xi(g)\rVert
_{\alpha,\,\lambda}  \\
& +C_{3}\lambda^{2\alpha-1}\Lambda_{\alpha}(h)\lVert
\xi(g)-\xi(h)\rVert_{\alpha,\,\lambda}.
\end{align*}
Hence the result.
\end{proof}

\medskip
Consider   Equation (\ref{f2}) on $\mathbb{R}^{d}$, for
$t\in\lbrack0,T]$, where $X_{0}\ $is a $d$-dimensional random variable, and
the coefficients $\sigma^{i,j},b^{i}:\Omega\times\left[  0,T\right]
\times\mathbb{R}^{d}\mathbb{\rightarrow}{}\mathbb{R}$ are measurable functions.

\begin{theorem}
\label{exist} Suppose that $X_{0}$ is an ${}\mathbb{R}^{d}$-valued random
variable, the coefficients $\sigma(t,x)$ and $b(t,x)$ satisfy assumptions (H1)
and (H2), where the constants might depend on $\omega$, with $\beta>1-H$,
$\delta>1/{H}-1$. Then if $\alpha\in\left(  1-H,\alpha_{0}\right)  $, then
there exists a unique stochastic process 
$$
X\in L^{0}\left(  \Omega
,\mathcal{F},\mathbb{P\,};W_{0}^{\alpha,\infty}(0,T;{}\mathbb{R}^{d})\right)
$$
 solution of the stochastic equation (\ref{f2}) and, moreover, for
$\mathbb{P}$-almost all $\omega\in\Omega$
\[
X\left(  \omega,\cdot\right)  =\left(  X^{i}\left(  \omega,\cdot\right)
\right)  _{d\times1}\in C^{1-\alpha}\left(  0,T;\mathbb{R}^{d}\right)  .
\]
\end{theorem}

Consider the particular case where $b=0$ and $\sigma$ is time independent,
that is,
\begin{equation}
X_{t}=X_{0}+\int_{0}^{t}\sigma(X_{s})dB_{s}^{H}. \label{3.1}%
\end{equation}
By the above theorem this equation has a unique solution provided $\sigma$ is
continuously differentiable, and $\sigma^{\prime}$ is bounded and H\"{o}lder
continuous of order $\delta>\frac{1}{H}-1$.

In \cite{HN} Hu and Nualart have established the following estimate. Choose
$\theta\in\left(  \frac{1}{2},H\right)  $. Then, the solution to Equation (\ref{3.1}) satisfies
\begin{equation}
\sup_{0\leq t\leq T}|X_{t}|\leq2^{1+kT\left(  \Vert\sigma^{\prime}%
\Vert_{\infty}\vee|\sigma(0)|\right)  \Vert B ^H\Vert_{{\theta}}^{1/{\theta}}%
}(|X_{0}|+1)\,. \label{eq1}%
\end{equation}
If $\sigma$ is bounded and $\Vert\sigma^{\prime}\Vert\not =0$ this estimate can be improved in the following way 
\begin{equation}
\sup_{0\leq t\leq T}|X_{t}|\leq|X_{0}|+k\Vert\sigma\Vert_{\infty}\left(
T^{\theta}\Vert B^H\Vert_{\theta}^{\frac{1}{\theta}}\vee T\Vert\sigma^{\prime}%
\Vert_{\infty}^{\frac{1-{\theta}}{{\theta}}}\Vert B^H\Vert_{{\theta}}%
^{\frac{1}{{\theta}}}\right)  . \label{eq2}%
\end{equation}
These estimates improve those obtained by Nualart and R\u{a}\c{s}canu in
\cite{NR} based on a suitable version of Gronwall's lemma. The estimates
(\ref{eq1}) and (\ref{eq2})  lead  to  the following integrability
properties for the solution of Equation (\ref{3.1}).

\begin{theorem}  \label{t1}
Consider the stochastic differential equation (\ref{3.1}), and assume that
$E(|X_{0}|^{p})<\infty$ for all $p\geq2$. If $\sigma^{\prime}$ is bounded and
H\"{o}lder continuous of order $\delta>\frac{1}{H}-1$, then
\begin{equation}
E\left(  \sup_{0\leq t\leq T}|X_{t}|^{p}\right)  <\infty
\end{equation}
for all $p\geq2$. If furthermore $\sigma$ is bounded and $E\left(
\exp(\lambda|X_{0}|^{\gamma})\right)  <\infty$ for any $\lambda>0$ and
$\gamma<2H$, then
\begin{equation}
E\left(  \exp\lambda\left(  \sup_{0\leq t\leq T}|X_{t}|^{\gamma}\right)
\right)  <\infty
\end{equation}
for any $\lambda>0$ and $\gamma<2H$.
\end{theorem}

In \cite{NS} Nualart and Saussereau have proved that the random variable
$X_{t}$ belongs locally to the space $\mathbb{D}^{\infty}$ if the function
$\sigma$ is infinitely differentiable and bounded together with all its
partial derivatives. As a consequence, they have derived the absolute
continuity of the law of $X_{t}$ for any $t>0$ assuming that the initial
condition is constant and the vector space spanned by $\{(\sigma^{i,}%
(x_{0}))_{1\leq i\leq d},1\leq j\leq m\}$ is $\mathbb{R}^{d}$.

Applying Theorem  \ref{t1} Hu and Nualart have proved in \cite{HN} that  if the function
$\sigma$ is infinitely differentiable and bounded together with all its
partial derivatives, then    for any $t\in [0,T]$ the random variable $X_{t}$ belongs to the space $\mathbb{D}^\infty$.  As a consequence, if the matrix $a(x)= \sigma \sigma^T(x)$ 
uniformly elliptic, then, for any $t>0$ the probability law of $X_{t}$ has an
$C^{\infty}$ density.  
In a recent paper, Baudoin and Coutin \cite{BC} have extended this result and
derived the regularity of the density under H\"{o}rmander hypoellipticity conditions.

\section{Flow of homeomorphisms}

\label{sec:main-results} Let $\pi=\left\{  0=t_{0}<t_{1}<\cdots<t_{n}%
=T\right\}  $ be the uniform partition of the interval $[0,T]$. That is
$t_{k}=\frac{kT}{n}$, $k=0,\ldots,n$. We denote by $B^{n,H\text{ }}$the
polygonal approximation of the fractional Brownian motion defined by%
\[
B_{t}^{n,H}=\sum_{k=0}^{n-1}\left(  B_{t_{k}}^{H}+\frac{n}{T}\left(
t-t_{k}\right)  \left(  B_{t_{k+1}}^{H}-B_{t_{k}}^{H}\right)  \right)
\mathbf{1}_{(t_{k},t_{k+1}]}(t).
\]

In order to get a precise rate for these approximations we will make use of
the following exact modulus of continuity of the fractional Brownian motion.
There exists a random variable $G$ such that almost surely for any
$s,t\in\lbrack0,T]$ we have%
\begin{equation}
\left|  B_{t}^{H}-B_{s}^{H}\right|  \leq G|t-s|^{H}\sqrt{\log\left(
|t-s|^{-1}\right)  }. \label{4.3b}%
\end{equation}
Fix $\theta<H$. We have the following result, which provides the rate of
convergence of these approximations in H\"{o}lder norm.

\begin{lemma}
\label{lema1} There exist a random variable $C_{T,\beta}$ such that
\begin{equation}
\Vert B^{H}-B^{n,H}\Vert_{{C}^{\theta}{(0,T;}\mathbb{R}^{m})}\leq C_{T,\beta
}n^{\theta-H}\sqrt{\log n}. \label{4.6}%
\end{equation}
\end{lemma}

\begin{proof}
To simplify the notation we will assume that $m=1$. Fix $0<s<t<T$ and assume
that $s\in\lbrack t_{l},t_{l+1}]$ and $t\in\lbrack t_{k},t_{k+1}]$. Let us
first estimate
\[
h_{1}(s,t)=\frac{1}{(t-s)^{\theta}}|B_{t}^{n,H}-B_{t}^{H}-(B_{s}^{n,H}%
-B_{s}^{H})|\,.
\]
If $t-s\geq\frac{T}{n}$, then using (\ref{4.3b}) we obtain
\begin{align*}
\left|  h_{1}(s,t)\right|   &  \leq T^{-\beta}n^{{\beta}}\left[  \left|
B_{t_{k}}^{H}-B_{t}^{H}+\frac{n}{T}\left(  t-t_{k}\right)  \left(  B_{t_{k+1}%
}^{H}-B_{t_{k}}^{H}\right)  \right|  \right. \\
&  \left.  +\left|  B_{t_{l}}^{H}-B_{s}^{H}+\frac{n}{T}\left(  s-t_{l}\right)
\left(  B_{t_{l+1}}^{H}-B_{t_{l}}^{H}\right)  \right|  \right] \\
&  \leq4GT^{-\theta+H}n^{-H+\theta}\sqrt{\log\left(  n/T\right)  }.
\end{align*}
If $t-s<\frac{T}{n}$, then there are two cases. Suppose first that
$s,t\in\lbrack t_{k},t_{k+1}]$. In this case, \ if $n$ is large enough we
obtain using (\ref{4.3b})%
\begin{align*}
\left|  h_{1}(s,t)\right|   &  \leq\frac{|B_{t}^{H}-B_{s}^{H}|}{(t-s)^{\theta
}}+\frac{n}{T}\frac{|B_{t_{k+1}}^{H}-B_{t_{k}}^{H}|}{(t-s)^{\theta}}(t-s)\\
&  \leq G|t-s|^{H-\theta}\sqrt{\log|t-s|^{-1}}+GT^{-1+H}\sqrt{\log
(n/T)}\ n^{1-H}(t-s)^{1-\theta}\\
&  \leq2GT^{-\theta+H}n^{-H+\theta}\sqrt{\log\left(  n/T\right)  }.
\end{align*}
On the other hand, if $s\in\lbrack t_{k-1},t_{k}]$ and $t\in\lbrack
t_{k},t_{k+1}]$ we have, again if $n$ is large enough
\begin{align*}
\left|  h_{1}(s,t)\right|   &  \leq\frac{1}{(t-s)^{\theta}}\left|  B_{t_{k}%
}^{H}-B_{t}^{H}+\frac{n}{T}\left(  t-t_{k}\right)  \left(  B_{t_{k+1}}%
^{H}-B_{t_{k}}^{H}\right)  \right. \\
&  \left.  -\left\{  B_{t_{k}}^{H}-B_{s}^{H}-\frac{n}{T}\left(  t_{k}%
-s\right)  \left(  B_{t_{k}}^{H}-B_{t_{k-1}}^{H}\right)  \right\}  \right| \\
&  \leq\frac{1}{(t-s)^{\theta}}\left[  |B_{t}^{H}-B_{s}^{H}|+\frac{n}%
{T}(t-s)\left(  |B_{t_{k}}^{H}-B_{t_{k-1}}^{H}|+|B_{t_{k+1}}^{H}-B_{t_{k}}%
^{H}|\right)  \right] \\
&  \leq\ \frac{G}{(t-s)^{\theta}}\left[  |t-s|^{H}\sqrt{\log|t-s|^{-1}%
}+2(t-s)\left(  \frac{n}{T}\right)  ^{H}\sqrt{\log\left(  n/T\right)  }\right]
\\
&  \leq3GT^{-\theta+H}n^{-H+\theta}\sqrt{\log\left(  n/T\right)  }.
\end{align*}
This proves (\ref{4.6}).
\end{proof}

\begin{corollary}
For any $\alpha\in(1-H,1/2)$, we have:
\[
\sup_{n}\Lambda_{\alpha}(B^{n,H})<+\infty\text{ and }\lim_{n\rightarrow
+\infty}\Lambda_{\alpha}(B^{n,H}-B^{H})=0.
\]
\end{corollary}

\begin{proof}
Choose $\eta>0$ in such a way that $1-\alpha+\eta<H$. According to
(\ref{eq:2}) and (\ref{eq:3}), we have
\[
\Lambda_{\alpha}(B^{n,H})\leq c_{\eta}\lVert B^{n,H}\rVert_{C^{1-\alpha+\eta
}(0,T;\mathbb{R}^{m})}%
\]
and%
\[
\Lambda_{\alpha}(B^{n,H}-B^{H})\leq c_{\eta}\lVert B^{n,H}-B^{H}%
\rVert_{C^{1-\alpha+\eta}(0,T;\mathbb{R}^{m})}.
\]
Then, Lemma \ref{lema1} implies that the sequence $B^{n,H}$ converges to
$B^{H\text{ }}$in the norm of \ $C^{1-\alpha+\eta}(0,T;\mathbb{R}^{m})$ which
yields the results.
\end{proof}

Consider for any $0\leq r\leq t\leq T$ and any natural number $n\geq1$ the
following equations
\begin{equation}
X_{rt}^{n}(x)=x+\int_{r}^{t}\sigma\left(  s,X_{rs}^{n}(x)\right)
dB^{n,H}(s)+\int_{r}^{t}b(s,X_{rs}^{n}(x))ds, \label{eq:4}%
\end{equation}
and
\begin{equation}
Y_{rt}^{n}(x)=x+\int_{r}^{t}\sigma\left(  s,Y_{st}^{n}(x)\right)
dB^{n,H}(s)+\int_{r}^{t}b(s,Y_{st}^{n}(x))ds. \label{eq:5}%
\end{equation}
We know from standard results on ordinary differential equations that for any
$n\geq1$,

\begin{enumerate}
\item \label{item:1} Equations (\ref{eq:4}) and (\ref{eq:5}) have a unique solution.

\item \label{item:2} For any $x\in{}\mathbb{R}^{d}$, for any $0\leq r\leq
\tau\leq t\leq T$, $X_{\tau t}^{n}(X_{r\tau}^{n}(x))=X_{rt}^{n}(x)$.

\item \label{item:3} For any $x\in\mathbb{R}^{d}$, for any $0\leq r\leq
\tau\leq t\leq T$, $Y_{r\tau}^{n}(Y_{\tau t}^{n}(x))=Y_{rt}^{n}(x)$.

\item \label{item:4} The maps $(x\mapsto X_{rt}^{n}(x))$ and $(x\mapsto
Y_{rt}^{n}(x))$ are $\mathbb{R}^{d}$-homeomorphisms inverse of each other:
\[
X_{rt}^{n}(Y_{rt}^{n}(x))=x\text{ and }Y_{rt}^{n}(X_{rt}^{n}(x))=x.
\]
\end{enumerate}

We are then in position to prove our main theorem:

\begin{theorem}
Assume that Hypothesis (H1) and (H2) hold. Then,
claims \ref{item:1}, \ref{item:2} \ref{item:3} and \ref{item:4} also hold for
\ the equations%
\begin{equation}
X_{rt}(x)=x+\int_{r}^{t}\sigma\left(  s,X_{rs}(x)\right)  dB^{H}(s)+\int
_{r}^{t}b(s,X_{rs}(x))ds, \label{g1}%
\end{equation}
and
\begin{equation}
Y_{rt}(x)=x+\int_{r}^{t}\sigma\left(  s,Y_{st}(x)\right)  dB^{H}(s)+\int
_{r}^{t}b(s,Y_{st}(x))ds. \label{g2}%
\end{equation}
\end{theorem}

\begin{proof}
Point \ref{item:1} is proved in \cite{NR}. As to the second claim,
proceed as follow:
\begin{align*}
X_{\tau t}^{n}(X_{r\tau}^{n}(x))-X_{\tau t}(X_{r\tau}(x))  &  =X_{\tau t}%
^{n}(X_{r\tau}^{n}(x))-X_{\tau t}^{n}(X_{r\tau}(x))\\
&  +(X_{\tau t}^{n}-X_{\tau t})(X_{r\tau}(x)).
\end{align*}
Fix $\varepsilon>0$ and $\alpha$ such that $1-H<\alpha<\frac{1}{2}$. Fix a
trajectory $\omega\in\Omega$. Choose $n_{0}$ so that $\Lambda_{\alpha}%
(B^{n,H}-B^{H})\leq\varepsilon$ for all $n\geq n_{0}$ and choose $\lambda$
such that $\lambda^{2\alpha-1}C_{2}\sup_{n}\Lambda_{\alpha}(B^{n,H}%
)\leq\frac{1}{2}$. Then, according to Theorem \ref{thm:continuite}, for any
$n\geq n_{0}$,
\[
\lVert X_{r\cdot}^{n}-X_{r\cdot}\rVert_{\alpha,\lambda}\leq\frac{C_{1}%
\lambda^{2\alpha-1}\lVert X_{r\cdot}\rVert_{\alpha,\,\lambda}}{1-C_{2}%
\lambda^{2\alpha-1}\ (1+\Lambda_{\alpha}(B^{n,H}))}\leq2C_{1}\lambda
^{2\alpha-1}\lVert X_{r\cdot}\rVert_{\alpha,\,\lambda}\varepsilon.
\]
Hence, for $n\geq n_{0}$,
\[
\left|  (X_{\tau t}^{n}-X_{\tau t})(X_{r\tau}(x))\right|  \leq c\varepsilon.
\]
The convergence of $X_{r\cdot}^{n}$ implies that there exists $R$ such that
for any $\tau\in\lbrack r,\,t]$ and for any $n\geq n_{0}$, $X_{r\tau}%
^{n}(x)\in B(0,R)$. Then, Theorem \ref{thm:continuite_condition_initiale}
implies that for $\lambda$ large enough%
\begin{align*}
  & \left|  X_{\tau t}^{n}(X_{r\tau}^{n}(x))-X_{\tau t}^{n}(X_{r\tau}(x))\right| \\
&   \quad \leq\left(  1-R\exp\left(  C\left(  1+\sup_{n}\Lambda_{\alpha}%
(B^{n,H})\right)  ^{\frac{1}{1-2\alpha}}\right)  \lambda^{2\alpha-1}\right)
^{-1}\\
&  \qquad  \times|X_{r\tau}^{n}(x)-X_{r\tau}(x)|\\
&   \quad \leq c|X_{r\tau}^{n}(x)-X_{r\tau}(x)|.
\end{align*}
We have thus proved that
\begin{align*}
0  &  =\lim_{n\rightarrow+\infty}X_{\tau t}^{n}(X_{r\tau}^{n}(x))-X_{\tau
t}(X_{r\tau}(x))\\
&  =\lim_{n\rightarrow+\infty}X_{rt}^{n}(x)-X_{\tau t}(X_{r\tau}(x))\\
&  =X_{rt}(x)-X_{\tau t}(X_{r\tau}(x)).
\end{align*}
Other points are handled similarly.
\end{proof}

\begin{acknowledgement}
This work was carried out during a stay of Laurent Decreusefond at Kansas
University, Lawrence KS. He would like to thank KU for warm hospitality and
generous support.
\end{acknowledgement}

%\bibliography{references}
%\bibliographystyle{amsplain}

\providecommand{\bysame}{\leavevmode\hbox to3em{\hrulefill}\thinspace}
\providecommand{\href}[2]{#2}

\end{document}